# Blowup of small data solutions for a quasilinear wave equation in two space dimensions

By Serge Alinhac

## Abstract

For the quasilinear wave equation

$$\partial_t^2 u - \Delta u = u_t u_{tt},$$

we analyze the long-time behavior of classical solutions with small (not rotationally invariant) data. We give a complete asymptotic expansion of the lifespan and describe the solution close to the blowup point. It turns out that this solution is a "blowup solution of cusp type," according to the terminology of the author [3].

## Résumé

Pour l'équation d'onde quasi-linéaire

$$\partial_t^2 u - \Delta u = u_t u_{tt},$$

nous analysons le comportement en grand temps des solutions classiques à données petites. Nous donnons un développement asymptotique complet du temps de vie et décrivons la solution près du point d'explosion. Cette solution est une "solution éclatée de type cusp," selon la terminologie de l'auteur [3].

## Introduction

We consider here the quasilinear equation in $\mathbf{R}^{2+1}$:

(0.1)
$$\partial_t^2 u - \Delta_x u = u_t u_{tt}$$

where

$$x_0 = t, \ x = (x_1, x_2), \ r = \sqrt{x_1^2 + x_2^2}, \ x_1 = r \cos \omega, \ x_2 = r \sin \omega.$$

We assume that the Cauchy data are $C^\infty$ and small,

$$u(x, 0) = \varepsilon u_1^0(x) + \varepsilon^2 u_2^0(x) + \ldots, \ u_t(x, 0) = \varepsilon u_1^1(x) + \varepsilon^2 u_2^1(x) + \ldots,$$

and supported in a fixed ball of radius $M$.

Our aim is to study the existence of smooth solutions to this problem, more precisely the lifespan $\bar{T}_\varepsilon$ of these solutions and the breakdown mechanism when these solutions stop being smooth.

This problem was introduced and extensively studied by John, for this and more general quasilinear wave equations, in space dimensions two or three (see his survey paper [9] and the references therein). Then lower bounds of the lifespan were obtained by Klainerman ([11], [12]), Hörmander ([7], [8]) and many other authors. Using some crude approximation by solutions of Burger's equation, Hörmander [7] has obtained in dimensions two and three explicit lower bounds for the lifespan. The result for equation (0.1) in dimension two is

$$(0.2) \qquad \liminf \varepsilon \bar{T}_\varepsilon^{1/2} \geq (\max \partial_\sigma^2 R^{(1)}(\sigma, \omega))^{-1} \equiv \bar{\tau}_0.$$

Here, the "first profile" $R^{(1)}$ is defined as

$$(0.3) \qquad R^{(1)}(\sigma, \omega) = \frac{1}{2\sqrt{2\pi}} \int_{s \geq \sigma} \frac{1}{\sqrt{s-\sigma}} [R(s, \omega, u_1^1) - \partial_s R(s, \omega, u_1^0)] ds,$$

where $R(s, \omega, v)$ denotes the Radon transform of the function $v$

$$R(s, \omega, v) = \int_{x\omega = s} v(x) dx.$$

Hörmander simply writes in his 1986 lectures on nonlinear hyperbolic equations [8]:

"Even if it is hard to doubt that (0.2) always gives the precise asymptotic lifespan of the solutions there is no proof except that of John [10] for the rotationally symmetric three-dimensional case."

In this paper, we prove Hörmander's conjecture that (0.2) indeed gives the correct asymptotic of the lifespan. In fact, our method of proof gives more than that : it provides a complete description of the solution close to the blowup point. It turns out that the solution is a "blowup solution of cusp type," according to the definitions of [3].

Finally, to formulate more precisely Hörmander's conjecture, let us introduce further useful notation and recall a previous result on upper bounds for the lifespan. Let $u_1$ be the solution of the linearized problem at 0:

$$\partial_t^2 u_1 - \Delta u_1 = 0, \quad u_1(x, 0) = u_1^0(x), \quad \partial_t u_1(x, 0) = u_1^1(x).$$

We have, for $r \to \infty, r - t \geq -C_0$, $R^{(1)}$ being the first profile defined by (0.3),

$$u_1 \sim \frac{R^{(1)}(r-t, \omega)}{r^{1/2}}.$$

Similarly, let us now define $u_2$ by

$$\partial_t^2 u_2 - \Delta u_2 - \partial_t u_1 \partial_t^2 u_1 = 0, \quad u_2(x, 0) = u_2^0(x), \quad \partial_t u_2(x, 0) = u_2^1(x).$$

We prove in [1] that, also for $r \to \infty, r - t \geq -C_0$,

$$u_2 - \frac{1}{2}(\partial_\sigma R^{(1)})^2 \sim \frac{R^{(2)}(r-t,\omega)}{r^{1/2}}$$

for a certain smooth $R^{(2)}$ that we call the "second profile." We assume that $\partial_\sigma^2 R^{(1)}$ has a unique positive quadratic maximum at a point $(\sigma_0, \omega_0)$, and then set

$$\bar{\tau}_0 = (\partial_\sigma^2 R^{(1)}(\sigma_0, \omega_0))^{-1},$$

$$\bar{\tau}_1 = -\bar{\tau}_0^2 \partial_\sigma^2 R^{(2)}(\sigma_0, \omega_0).$$

The result of [2] (which is also valid for general quasilinear wave equations) is the following.

ASYMPTOTIC THEOREM (see [2]).   *Under the above nondegeneracy assumption on the initial data, there exists a function $\bar{T}_\varepsilon^a$ with the following properties*:
i) *For all $N$, $\bar{T}_\varepsilon \geq \bar{T}_\varepsilon^a - \varepsilon^N$ for $0 < \varepsilon \leq \varepsilon_N$,*
ii) *For some $C > 0$ and $(C\varepsilon^2)^{-1} \leq t \leq \bar{T}_\varepsilon^a - \varepsilon^N$,*

$$\frac{1}{C}\frac{1}{\bar{T}_\varepsilon^a - t} \leq |\nabla^2 u(.,t)|_{L^\infty} \leq C\frac{1}{\bar{T}_\varepsilon^a - t}.$$

*The function $\bar{T}_\varepsilon^a$ is of the form*

$$\bar{T}_\varepsilon^a = \varepsilon^{-2}(\bar{\tau}_\varepsilon^a)^2(\varepsilon, \varepsilon^2 ln\varepsilon),$$

*where $\bar{\tau}_\varepsilon^a$ is a smooth function satisfying*

$$\bar{\tau}_\varepsilon^a = \bar{\tau}_0 + \varepsilon\bar{\tau}_1 + O(\varepsilon^2 ln\varepsilon).$$

Thus, for numerical purposes, the *asymptotic lifespan* $\bar{T}_\varepsilon^a$ looks like the true lifespan $\bar{T}_\varepsilon$; this feature would certainly make numerical experiments, designed to test whether or not the solution actually blows up at time $\bar{T}_\varepsilon^a$, very hard to realize.

We prove in this work that, for equation (0.1), one has in fact $\bar{T}_\varepsilon \sim \bar{T}_\varepsilon^a$.

## I. Results and method of proof

1. Throughout this paper, we make the following nondegeneracy assumption on the initial data.

**(ND) The function $\partial_\sigma^2 R^{(1)}(\sigma, \omega)$ has a unique positive quadratic maximum at a point $(\sigma_0, \omega_0)$.**

Recall that the first profile $R^{(1)}$ was defined in (0.3).

For equation (0.1) with small data satisfying (ND), we have the following theorem.

LIFESPAN THEOREM 1.1.1.    *The lifespan $\bar{T}_\varepsilon$ of the solution $u$ of (0.1) satisfies*

$$\bar{\tau}_\varepsilon \equiv \varepsilon(\bar{T}_\varepsilon)^{1/2} = \bar{\tau}_0 + \varepsilon\bar{\tau}_1 + O(\varepsilon^2 ln\varepsilon). \tag{1.1.1}$$

*Moreover, for $t \geq \tau_0^2 \varepsilon^{-2}$ ($0 < \tau_0 < \bar{\tau}_0$) and $\varepsilon$ small,*
  i) *The solution $u$ is of class $C^1$ and $|u|_{C^1} \leq C\varepsilon^2$;*
  ii) *There is a point $M_\varepsilon = (m_\varepsilon, \bar{T}_\varepsilon)$ such that, away from $M_\varepsilon$, the solution $u$ is of class $C^2$ with $|u|_{C^2} \leq C\varepsilon^2$ there;*
  iii) *The solution satisfies*

$$|\nabla^2 u(.,t)|_{L^\infty} \leq \frac{C}{\bar{T}_\varepsilon - t}, \tag{1.1.2}$$

$$|\partial_t^2 u(.,t)|_{L^\infty} \geq \frac{1}{C}\frac{1}{\bar{T}_\varepsilon - t}. \tag{1.1.3}$$

We give here only the approximation (1.1.1) for simplicity. In fact, it is easily seen that the lifespan $\bar{T}_\varepsilon$ and the location of the blowup point $M_\varepsilon$ can be computed to any order (for small enough $\varepsilon$) by the implicit function arguments of [2]. In particular, $\bar{T}_\varepsilon \sim \bar{T}_\varepsilon^a$ in the sense of asymptotic series.

The inequalities (1.1.2), (1.1.3) give a rough idea of how the second order derivatives of the solution blow up. A much better description of the solution close to $M_\varepsilon$ can be obtained from the following theorem.

GEOMETRIC BLOWUP THEOREM 1.1.2.    *There exist a point $\tilde{M}_\varepsilon = (\tilde{m}_\varepsilon, \bar{\tau}_\varepsilon)$, a neighbourhood $V$ of $\tilde{M}_\varepsilon$ in $\{(s, \omega, \tau), s \in R, \omega \in S^1, \tau \leq \bar{\tau}_\varepsilon\}$ and functions $\phi, \tilde{G}, \tilde{v} \in C^3(V)$ with the following properties:*
  i) *The function $\phi$ satisfies in $V$ the condition*

$$\phi_s \geq 0, \phi_s(s, \omega, \tau) = 0 \Leftrightarrow (s, \omega, \tau) = \tilde{M}_\varepsilon, \tag{H}$$
$$\phi_{s\tau}(\tilde{M}_\varepsilon) < 0, \nabla_{s,\omega}(\phi_s)(\tilde{M}_\varepsilon) = 0, \nabla_{s,\omega}^2(\phi_s)(\tilde{M}_\varepsilon) >> 0.$$

  ii) *$\partial_s \tilde{G} = \phi_s \tilde{v}$ and $\partial_s \tilde{v}(\tilde{M}_\varepsilon) \neq 0$. If we define the map*

$$\Phi(s, \omega, \tau) = (\sigma = \phi(s, \omega, \tau), \omega, \tau)$$

*and set $\Phi(\tilde{M}_\varepsilon) = (|x_\varepsilon| - \bar{T}_\varepsilon, x_\varepsilon|x_\varepsilon|^{-1}, \bar{\tau}_\varepsilon) \equiv \bar{M}_\varepsilon$, condition (H) allows us to define near $\bar{M}_\varepsilon$ a function $G$ by*

$$G(\Phi) = \tilde{G}. \tag{1.1.4}$$

*Then, close to $M_\varepsilon = (x_\varepsilon, \bar{T}_\varepsilon)$, the solution $u$ satisfies*

$$u(x,t) = \frac{\varepsilon}{r^{1/2}}G(r - t, \omega, \varepsilon t^{1/2}). \tag{1.1.5}$$

*Finally, the functions $\phi, \tilde{G}$ and $\tilde{v}$ are of class $C^k$ for $\varepsilon \leq \varepsilon_k$.*

In this theorem, we see that the singularities of $u$ come only from the singularities of $G$ ; these in turn arise from the fact that the mapping $\Phi$ is not invertible at the point $\tilde{M}_\varepsilon$. More precisely, condition (H) implies that the singularity of $\Phi$ is a *cusp singularity*. Thus, describing the behavior of the derivatives of $u$ near $M_\varepsilon$ is just a local geometric problem. This is the reason why we call this behavior of $u$ "geometric blowup" (see [3] or [5] for details).

2. Let us explain now the method of proof of Theorems 1.1.1 and 1.1.2. The idea is to construct a piece of blowup solution to (0.1) in a strip

$$-C_0 \leq r - t \leq M, \ \tau_0^2 \varepsilon^{-2} \leq t \leq \bar{T}_\varepsilon, \ 0 < \tau_0 < \bar{\tau}_0$$

close to the boundary of the light cone. This gives an upper bound for the lifespan, which turns out to be the correct one. Of course, this is not surprising, because the first blowup of the solution is believed to take place in such a strip, and not far inside the light cone.

The proof is thus devoted to this construction, which is done in four steps, handled respectively in parts II, III, IV and V of the present paper.

*Step* 1: *Asymptotic analysis, normalization of variables and reduction to a local problem.* We choose a number $0 < \tau_0 < \bar{\tau}_0$ and use here asymptotic information on the behavior of $u$ for $r - t \geq -C_0$ and $\varepsilon t^{1/2}$ close to $\tau_0$. Thus, we are far away from any possible blowup at this stage, because of (0.2). According to [1], the solution in this domain behaves like a smooth function (depending smoothly also on $\varepsilon$ and $\varepsilon^2 ln\varepsilon$) of the variables

$$\sigma = r - t, \ \omega, \ \tau = \varepsilon t^{1/2}.$$

Thus we set

$$u(x,t) = \frac{\varepsilon}{r^{1/2}} G(\sigma, \omega, \tau).$$

Writing equation (0.1) for $G$ in these new variables, we are left with solving a local problem for $G$ in a domain

$$-C_0 \leq \sigma \leq M, \ \tau_0 \leq \tau \leq \bar{\tau}_\varepsilon,$$

where $\bar{\tau}_\varepsilon = \varepsilon \bar{T}_\varepsilon^{1/2}$ is still unknown. At this stage, we have a *free boundary problem*, the upper boundary of the domain being determined by the first blowup time.

*Step* 2: *Blowup of the problem.* To solve the free boundary problem of Step 1, we introduce a singular (still unknown) change of variables

$$\Phi : (s, \omega, \tau) \mapsto (\sigma = \phi(s, \omega, \tau), \omega, \tau), \ \phi(s, \omega, \tau_0) = s.$$

The idea is to obtain $G$ in the form

$$G(\Phi) = \tilde{G}$$

for *smooth* functions $\phi$ and $\tilde{G}$, and arrange at the same time to have $\phi_s$ vanish at one point $\tilde{M}_\varepsilon = (\tilde{m}_\varepsilon, \bar{\tau}_\varepsilon)$ of the upper boundary of the domain. Thus, we will have

$$\tilde{G}_s = G_\sigma \phi_s,$$

and the technical condition ii) of Theorem 1.1.2 gives in fact

$$G_\sigma(\Phi) = \tilde{v};$$

hence

$$G_{\sigma\sigma}(\Phi) = \tilde{v}_s/\phi_s.$$

We see that $u, \nabla u$ will remain continuous and that $\nabla^2 u$ will blow up at some point, in accordance with the expected behavior of $u$.

The nonlinear system on $\phi$ and $\tilde{G}$ corresponding to (0.1) is called the blowup system.

Instead of looking for a *singular* solution of the normalized original equation as in Step 1, we are now looking for a *smooth* solution of the blowup system ; however, we cannot just solve for $\tau$ close to $\tau_0$ : we have to reach out to attain a point where $\phi_s = 0$.

Finally, introducing an unknown real parameter (corresponding to the height of the domain), we can reduce the free boundary problem of Step 1 to a problem in a fixed domain.

*Step* 3: *Existence and tame estimates for a linear Goursat problem.* Linearization of the problem obtained in Step 2 leads to a third order Goursat problem. In fact, it is the special structure of (0.1) which makes it possible to reduce the full blowup system on $\phi$ and $\tilde{G}$ to a *scalar* equation on $\phi$. The (unknown) point where $\phi_s$ vanishes is a degeneracy point for this equation. Energy estimates can then be obtained using an appropriate multiplier. We prove in this step existence of solutions and tame estimates, which allow us to solve the nonlinear problem by a Nash-Moser method.

*Step* 4: *Back to the solution $u$.* Having $\tilde{G}$ and $\phi$, we deduce $G$ and thus obtain a piece of solution $\tilde{u}$ of (0.1) with the desired properties. It remains to see that $\tilde{u} = u$ where $\tilde{u}$ is defined, and that $u$ does not blow up anywhere else.

## II. Step 1: Asymptotic analysis, normalization of variables and reduction to a local problem

1. The asymptotic analysis of (0.1) was carried out in [1]. Fix

$$0 < \tau_1 < \tau_0 < \tau_2 < \bar{\tau}_0.$$

Introducing the variables

$$\sigma = r - t, \ \omega, \ \tau = \varepsilon t^{1/2}$$

as before, we only need here the behavior of the solution in the region

$$\{\tau_1 \leq \tau \leq \tau_2, \ r - t \geq -C_0\},$$

that is, far away from any possible blowup. The result of [1] is that if we set

$$u(x,t) = \frac{\varepsilon}{r^{1/2}} G(\sigma, \omega, \tau),$$

the function $G$ is bounded in $C^k$ (independently of $\varepsilon$) for $\varepsilon \leq \varepsilon_k$ ($\varepsilon_k$ depends of course on $C_0, \tau_1$ and $\tau_2$). For $\varepsilon = 0$, the function $G$ reduces to the function, abusively denoted by $R^{(1)}(\sigma, \omega, \tau)$, solution of the Cauchy problem

$$(2.1.1) \qquad \partial_\tau G - \frac{1}{2}(\partial_\sigma G)^2 = 0, \ G(\sigma, \omega, 0) = R^{(1)}(\sigma, \omega).$$

According to a simple computation, the function $G$ satisfies an equation of the form

$$(2.1.2) \qquad -\partial^2_{\sigma\tau} G + (\partial_\sigma G)(\partial^2_\sigma G) + \varepsilon^2 E(\sigma, \omega, \tau, G, \nabla G, \nabla^2 G) = 0,$$

where $E$ is a smooth function, linear in $\nabla^2 G$, which we need not know explicitly.

2. To prepare for Step 2, it is important to see that if we set $w = u_t$ and take the t-derivative of the equation, we obtain the *conservative* nonlinear equation

$$(2.2.1) \qquad P(w) = \partial_t^2 w - \Delta w - \frac{1}{2}\partial_t^2(w^2) = 0.$$

Note that, with $w = \frac{\varepsilon}{r^{1/2}} F$,

$$(2.2.2) \qquad F = \mathcal{L}_1 G, \ \mathcal{L}_1 = -\partial_\sigma + \frac{\varepsilon^2}{2\tau}\partial_\tau.$$

We need the expression of $P(w)$ in the variables $\sigma, \omega, \tau$.

LEMMA II.2.  *There exists the identity*

$$(2.2.3) \quad \frac{r}{\varepsilon^2} P(w) = -FF_{\sigma\sigma} - \frac{R^{1/2} - \varepsilon^2 F}{\tau}[F_{\sigma\tau} - \frac{\varepsilon^2}{4\tau}F_{\tau\tau}]$$
$$- \varepsilon^2 R^{-3/2} F_{\omega\omega} - (F_\sigma - \frac{\varepsilon^2}{2\tau}F_\tau)^2 + \varepsilon^2 h \nabla F + \varepsilon^2 h_0 \equiv \tilde{P}(F)$$

where $h$ and $h_0$ are smooth functions of $(\omega, R, \tau, F)$, and $R = \tau^2 + \varepsilon^2 \sigma$.

We want to solve $\tilde{P}(F) = 0$ in a (still unknown ) domain

$$-A_0 \leq \sigma, \ \omega \in S^1, \ \tau_0 \leq \tau \leq \bar{\tau}_\varepsilon,$$

with two trace conditions on $\{\tau = \tau_0\}$ corresponding to that for $u$ and $F$ supported in $\{\sigma \leq M\}$ ($A_0$ is big enough).

### III. Step 2: Blowup of the problem and reduction to a Goursat problem on a fixed domain

1. *Formal blowup.* We set, with an unknown $\phi$,

(3.1.1) $$G(\Phi) = \tilde{G}, \quad F(\Phi) = v, \quad \Phi(s, \omega, \tau) = (\phi(s, \omega, \tau), \omega, \tau).$$

We have then, with $y = \omega$ or $\tau$,

(3.1.2) $$(\partial_\sigma G)(\Phi) = \phi_s^{-1} \partial_s \tilde{G}, \quad (\partial_y G)(\Phi) = \partial_y \tilde{G} - \left(\frac{\phi_y}{\phi_s}\right) \partial_s \tilde{G},$$

and in particular

(3.1.3) $$\mathcal{L}_1 G(\Phi) = \phi_s^{-1} \bar{\mathcal{L}}_1 \tilde{G}, \quad \bar{\mathcal{L}}_1 = -\left(1 + \frac{\varepsilon^2 \phi_\tau}{2\tau}\right)\partial_s + \frac{\varepsilon^2 \phi_s}{2\tau}\partial_\tau.$$

For second order derivatives of $G$, we find an expression of the form

(3.1.4) $$(\nabla^2 G)(\Phi) = \frac{\phi_{ss}}{\phi_s^3} A - \frac{\partial_s A}{\phi_s^2} + \frac{B}{\phi_s},$$

where $A$ and $B$ are smooth.

Let us explain now heuristically how we establish the blowup system. Our aim is to substitute the expressions (3.1.2) and (3.1.4) into the equation (2.1.2) for $G$ and take the coefficients of the various powers of $\phi_s^{-1}$ to be zero. Of course, if we do this in a straightforward manner, we will obtain too many equations on $\tilde{G}$ and $\phi$. Another possibility is to introduce an auxiliary (smooth) function $\tilde{v}$ and force the relation

(3.1.5) $$\partial_s \tilde{G} = \phi_s \tilde{v}.$$

We see then from (3.1.2) that $\nabla G$ is smooth and $\nabla^2 G$ is of the form $\frac{A}{\phi_s} + B$ (with $A, B$ smooth); equating to zero the coefficients of $1$ and of $\phi_s^{-1}$ in the equation for $G$ yields then two equations, which give, along with (3.1.5), a $(3 \times 3)$-system on $\tilde{G}, \tilde{v}, \phi$. Here, we take advantage of formula (3.1.4) and of the conservative character of equation (2.2.1) to get a $(2 \times 2)$-system on $v, \phi$, as indicated in the following lemma.

LEMMA III.1.  *Since the functions $v$ and $F$ are related by (3.1.1),*

$$\tilde{P}(F)(\Phi) = \frac{1}{\phi_s^3} \phi_{ss} v_s \mathcal{T}_0 + \frac{1}{\phi_s^2} \mathcal{T}_1 + \frac{1}{\phi_s} \mathcal{T}_2 + \mathcal{T}_3,$$

*where*

(3.1.6) $$\mathcal{T}_0 = qv - \frac{R^{1/2}}{\tau}\phi_\tau - \varepsilon^2 \frac{R^{1/2}}{4\tau^2}\phi_\tau^2 + \frac{\varepsilon^2}{R^{3/2}}\phi_\omega^2,$$

(3.1.7) $$\mathcal{T}_1 = -\partial_s(v_s \mathcal{T}_0),$$

(3.1.8) $$\mathcal{T}_2 = Z\partial_s v - \varepsilon^2 v_s N\phi + \varepsilon^2 v_s h_2(\omega, \tau, v, v_\omega, v_\tau, \phi, \phi_\omega, \phi_\tau),$$

(3.1.9) $$\mathcal{T}_3 = \varepsilon^2 Nv + \varepsilon^2 h_3(\omega, \tau, v, v_\omega, v_\tau, \phi),$$

$h_2$ and $h_3$ being smooth functions and $R = \tau^2 + \varepsilon^2 \phi$. Moreover,

$$(3.1.10) \quad q = 1 + \frac{\varepsilon^2}{\tau}\phi_\tau + \frac{\varepsilon^4}{4\tau^2}\phi_\tau^2, \quad Z = \delta_1 \partial_\tau + \varepsilon^2 \delta_2 \partial_\omega$$

with

$$(3.1.11) \quad \delta_1 = -\frac{1}{\tau}\{R^{1/2} - \varepsilon^2 v + \varepsilon^2 \frac{R^{1/2} - \varepsilon^2 v}{2\tau}\phi_\tau\}, \quad \delta_2 = 2R^{-3/2}\phi_\omega.$$

$$(3.1.12) \quad N = \frac{R^{1/2} - \varepsilon^2 v}{4\tau^2}\partial_\tau^2 - R^{-3/2}\partial_\omega^2 \equiv N^{(1)}\partial_\tau^2 + 2\varepsilon^2 N^{(2)}\partial_{\tau\omega}^2 + N^{(3)}\partial_\omega^2.$$

We note the three following important facts:

$$(3.1.13) \quad q \neq 0,$$
$$(3.1.14) \quad \delta_1 = -1 + O(\varepsilon^2),$$
$$(3.1.15) \quad N^{(1)} = \frac{1}{4\tau} + O(\varepsilon^2) > 0, \quad N^{(3)} = -\frac{1}{\tau^3} + O(\varepsilon^2) < 0.$$

The fact that $N^{(2)}$ is actually zero does not play a role in the subsequent computations, so that it is more natural to keep it.

In order to solve the equation $\tilde{P}(F) = 0$, we now take $v$ and $\phi$ to solve the blowup system

$$(3.1.16) \quad \mathcal{T}_0 = 0, \quad \mathcal{T}_2 + \phi_s \mathcal{T}_3 = 0.$$

2. *Reduction to a free boundary Goursat problem.* In this section, we are going to reduce the blowup system (3.1.16) to a scalar problem on $\phi$, with boundary conditions given on characteristic boundaries of the (still unknown) domain.

2.1. *A local solution of the blowup system.* From the implicit function theorem, we can write equation $\mathcal{T}_0 = 0$ in the form

$$\phi_\tau = E(\omega, \tau, \phi, \phi_\omega, v),$$

with

$$(3.2.1) \quad E(\omega, \tau, \phi, 0, 0) = 0,$$

and, for $\varepsilon = 0$,

$$E = v.$$

The function $F$ being in fact known and smooth in a small strip

$$S_1 = \{\tau_0 \leq \tau \leq \tau_0 + \eta, \eta > 0\},$$

we can solve, for $\eta$ small enough, the Cauchy problem

$$\phi_\tau = E(\omega, \tau, \phi, \phi_\omega, F(\phi, \omega, \tau)), \quad \phi(s, \omega, \tau_0) = s$$

in this strip. Setting then $v = F(\phi, \omega, \tau)$, we obtain a smooth particular solution $(\bar{v}, \bar{\phi})$ of (3.1.16). Note that, thanks to (3.2.1), $\bar{v}$ and $\bar{\phi} - s$ are smooth and flat on $\{s = M\}$.

2.2. *Straightening out of a characteristic surface.* Consider the "nearly horizontal" surface $\Sigma$ through $\{\tau = \tau_0, s = M\}$ which is characteristic for the operator $Z\partial_s + \varepsilon^2 \bar{\phi}_s N$, the coefficients of $Z$ and $N$ being computed on $(\bar{v}, \bar{\phi})$. The surface $\Sigma$ is defined by an equation

$$\tau = \psi(s, \omega) + \tau_0,$$

where $\psi$ is the solution of the Cauchy problem
(3.2.2)
$$(-\delta_1 + \varepsilon^2 \delta_2 \psi_\omega)\psi_s + \varepsilon^2 \bar{\phi}_s(N^{(1)} - 2\varepsilon^2 N^{(2)}\psi_\omega + N^{(3)}\psi_\omega^2) = 0, \quad \psi(M, \omega) = 0.$$

Equation (3.2.2) has, for small $\varepsilon$, a smooth solution in the appropriate domain. This solution is $O(\varepsilon^2)$ and decreasing in $s$.

We now perform the change of variables

$$(3.2.3) \quad X = s, \quad Y = \omega, \quad T = (1 - \chi(\frac{\tau - \tau_0}{\eta}))(\tau - \tau_0) + (\tau - \tau_0 - \psi)\chi(\frac{\tau - \tau_0}{\eta}),$$

where $\chi \in C^\infty$, $\chi(t) = 1$ for $t \leq 1/2$, $\chi(t) = 0$ for $t \geq 1$. Note that this change reduces to $T = \tau - \tau_0$ away from a neighbourhood of $\{\tau = \tau_0\}$. The (still unknown) domain

$$D_\psi = \{-A_0 \leq s \leq M, \omega \in S^1, \tau_0 + \psi \leq \tau \leq \bar{\tau}_\varepsilon\}$$

is taken by this change into

$$\tilde{D} = \{-A_0 \leq X \leq M, Y \in S^1, 0 \leq T \leq \bar{T} = \bar{\tau}_\varepsilon - \tau_0\}.$$

With a slight abuse of notation, we will again denote by $(\bar{v}, \bar{\phi})$ the local solution of (3.1.14) transformed by (3.2.3); this solution exists now in a small strip $\{0 \leq T \leq \eta_1\}$ of $\tilde{D}$.

2.3. *Reduction to an equation on $\phi$.* The equation $\mathcal{T}_0 = 0$ allows us to express $v$ in terms of $\phi$ in the form

$$(3.2.4) \quad v = V(\omega, \tau, \phi, \phi_\omega, \phi_\tau).$$

Replacing $v$ by $V$ in (3.1.16), we obtain a third order equation on $\phi$, according to Lemma III.1. The change of variables (3.2.3) gives

$$\partial_s = \partial_X + T_s \partial_T \equiv S, \quad \partial_\omega = \partial_Y + T_\omega \partial_T, \quad \partial_\tau = T_\tau \partial_T,$$

where

$$T_s = O(\varepsilon^2), \quad T_\omega = O(\varepsilon^2), \quad T_\tau = 1 + O(\varepsilon^2)$$

are known functions. Hence the equation on $\phi$ becomes, in the new variables,

(3.2.5) $\quad \mathcal{L}(\phi) \equiv (\tilde{Z}S)\tilde{V} + \varepsilon^2(S\phi)\tilde{N}\tilde{V} - \varepsilon^2(S\tilde{V})\tilde{N}\phi + \varepsilon^2(S\tilde{V})\tilde{h}_2 + \varepsilon^2(S\phi)\tilde{h}_3 = 0,$

where $\tilde{Z}$, $\tilde{N}$, $\tilde{V}$ etc. correspond to $Z$, $N$, $V$ etc., transformed by (3.2.3). We note then

(3.2.6) $$\tilde{Z} = \tilde{\delta}_1 \partial_T + \varepsilon^2 \tilde{\delta}_2 \partial_Y,$$

(3.2.7) $$\tilde{N} = \tilde{N}^{(1)} \partial_T^2 + 2\varepsilon^2 \tilde{N}^{(2)} \partial_{YT}^2 + \tilde{N}^{(3)} \partial_Y^2.$$

Our goal is now to solve $\mathcal{L}(\phi) = 0$ in $\tilde{D}$ with the boundary conditions

(3.2.8) $\quad \phi(X, Y, 0) = \bar{\phi}(X, Y, 0), \quad \partial_T \phi(X, Y, 0) = \partial_T \bar{\phi}(X, Y, 0),$

and $\phi - X$ is flat on $\{X = M\}$.

**2.4. Construction of an approximate solution in the large.** Note that for $\varepsilon = 0$, the change (3.2.3) reduces to the translation $T = \tau - \tau_0$, while the blowup system (3.1.16) is

$$v = \phi_\tau, \quad \partial_{\tau s}^2 v = 0.$$

The initial conditions for this system are

$$\phi(X, Y, 0) = X, \quad \partial_T \phi(X, Y, 0) = -\partial_\sigma R^{(1)}(X, Y, \tau_0).$$

Hence the value $\bar{\phi}_0$ of $\bar{\phi}$ for $\varepsilon = 0$ is

$$\bar{\phi}_0(X, Y, T) = X - T R_\sigma^{(1)}(X, Y, \tau_0).$$

To obtain an approximate solution valid also for large values of $T$, we just glue together the local true solution $\bar{\phi}$ to $\bar{\phi}_0$:

$$\bar{\phi}^{(0)}(X, Y, T) = \chi(\frac{T}{\eta_1})\bar{\phi}(X, Y, T) + (1 - \chi(\frac{T}{\eta_1}))\bar{\phi}_0(X, Y, T).$$

We have then

$$\mathcal{L}(\bar{\phi}^{(0)}) = \bar{f}^{(0)},$$

where $\bar{f}^{(0)}$ is smooth, flat on $\{X = M\}$, zero near $\{T = 0\}$, and zero for $\varepsilon = 0$.

**2.5. The condition (H).** Let us consider more closely the vanishing of $\phi_X$ in $\tilde{D}$. On one hand, $\phi_X$ has to vanish somewhere, otherwise the corresponding $F$ and $u$ would not have any singularities. On the other hand, as will be clear from the linear analysis of Chapter IV, the linearized problem corresponding to $\mathcal{L}(\phi) = 0$ seems to become unstable for $\phi_X < 0$. Hence we are forced to consider the situation where $\phi_X$ vanishes only on the upper boundary of $\tilde{D}$.

In a way completely analogous to what we have done in [4], we expect $\phi$ to satisfy, for some point $\tilde{M} = (\tilde{m}, \bar{T})$, the condition

(H) $$\phi_X \geq 0, \ \phi_X(X, Y, T) = 0 \Leftrightarrow (X, Y, T) = \tilde{M},$$
$$\phi_{XT}^2(\tilde{M}) < 0, \ \nabla_{X,Y}(\phi_X)(\tilde{M}) = 0, \ \nabla_{X,Y}^2(\phi_X)(\tilde{M}) \gg 0.$$

Let us show that the approximate solution $\bar{\phi}^{(0)}$ from 2.4 satisfies this condition (H) at time

(3.2.9) $$\bar{T} = T_0 = (\max \partial_X^2 R^{(1)}(X, Y, \tau_0))^{-1}.$$

Thanks to the nondegeneracy assumption (ND), the function

$$\partial_\sigma(-\partial_\sigma R^{(1)}(\sigma, \omega))$$

has a quadratic minimum at $(\sigma_0, \omega_0)$. On the other hand, the function $-\partial_\sigma R^{(1)}(\sigma, \omega, \tau)$ is a solution of Burger's equation: at time $\tau_0$, its $\sigma$ derivative also has a quadratic minimum at the corresponding point, image of $(\sigma_0, \omega_0)$ by the characteristic flow. In addition, $T_0 = \bar{\tau}_0 - \tau_0$. Finally, $\partial_X \bar{\phi} > 0$ close to $\{T = 0\}$.

3. *Reduction to a Goursat problem on a fixed domain and condition* (H).

3.1. *Reduction to a fixed domain.* Recall that we want to solve the equation $\mathcal{L}(\phi) = 0$ in a domain such that $\phi$ satisfies the condition (H) for a point located on the upper boundary. The approximate solution $\bar{\phi}^{(0)}$, starting point of some approximation process, satisfies this condition for a domain of height $T_0$, according to 2.4, 2.5. Unfortunately, in the successive approximation process, further modifications of $\bar{\phi}^{(0)}$ will yield functions not satisfying (H) anymore. We are thus forced, at each step of the process, to adjust the domain to have the new $\phi$ satisfy condition (H).

To achieve this, we introduce a real parameter $\lambda$ close to zero, and perform the change of variables

(3.3.1) $$X = x, \ Y = y, \ T \equiv T(t, \lambda) = T_0(t + \lambda t(1 - \chi_1(t))),$$

where $\chi_1$ is 1 near 0 and 0 near 1, and $T_0$ is defined as in (3.2.9). Of course, one should not confuse these variables with the original variables ! We will from now on work on a fixed domain

$$D_0 = \{-A_0 \leq x \leq M, \ y \in S^1, \ 0 \leq t \leq 1\}.$$

We denote the transformed equation by

(3.3.2) $$\mathcal{L}(\lambda, \phi) = 0,$$

the transformed approximate solution for $\lambda = \lambda^{(0)} = 0$ by

(3.3.3) $$\phi^{(0)}(x, y, t) = \bar{\phi}^{(0)}(x, y, T(t, 0)) = \bar{\phi}^{(0)}(x, y, T_0 t),$$

and set

(3.3.4) $$\mathcal{L}(\lambda^{(0)}, \phi^{(0)}) = f^{(0)} \equiv \bar{f}^{(0)}(x, y, T_0 t).$$

We note that $\phi^{(0)}$ satisfies (H) in $D_0$ for a certain point $\tilde{M}_0 = (\tilde{m}_0, 1)$.

3.2. *Structure of the linearized operator.* The linearized operator of $\mathcal{L}$ at the point $(\lambda, \phi)$ is denoted by

(3.3.5) $$\mathcal{L}'(\lambda, \phi)(\dot{\lambda}, \dot{\phi}) = \partial_\lambda \mathcal{L}(\lambda, \phi)\dot{\lambda} + \partial_\phi \mathcal{L}(\lambda, \phi)\dot{\phi}.$$

Because $\mathcal{L}(\lambda, \phi)$ comes from $\mathcal{L}(\phi)$ by (3.2.1), we have the following lemma.

LEMMA III.3.1. *If $\mathcal{L}(\lambda, \phi) = f$, then*

(3.3.6) $$\partial_\lambda \mathcal{L}(\lambda, \phi) + \partial_\phi \mathcal{L}(\lambda, \phi)(\partial_t \phi \frac{\partial_\lambda T}{\partial_t T}) = \partial_t f \frac{\partial_\lambda T}{\partial_t T}.$$

For the time being, it is not necessary to make an explicit computation of $\partial_\phi \mathcal{L}$. Note only that, if we have at some stage $\mathcal{L}(\lambda, \phi) = f$ (for a small $f$), to solve

$$\mathcal{L}'(\lambda, \phi)(\dot{\lambda}, \dot{\phi}) = \dot{f}$$

approximately it is enough to solve

$$\partial_\phi \mathcal{L}(\lambda, \phi)\dot{\Psi} = \dot{f}$$

and then to take $(\dot{\lambda}, \dot{\phi})$, verifying

$$\dot{\phi} - \dot{\lambda}\phi_t \frac{\partial_\lambda T}{\partial_t T} = \dot{\Psi}.$$

In fact, we get with this choice

$$\mathcal{L}'(\lambda, \phi)(\dot{\lambda}, \dot{\phi}) = \dot{f} - \dot{\lambda}\partial_t f \frac{\partial_\lambda T}{\partial_t T}.$$

The additionnal term contains a product of small $\dot{\lambda}$ by small $\partial_t f$, which is negligible as a quadratic error. Having determined $\dot{\Psi}$, we see that we still have an additional degree of freedom to choose $\dot{\phi}$: we will take advantage of this to arrange for $\phi + \dot{\phi}$ to satisfy (H).

3.3. *The fundamental lemma.* We follow here exactly the same idea as in [4].

LEMMA III.3.2. *Assume that $\phi - \phi^{(0)}$ and $\psi$ are small enough in $C^4(D_0)$. Then*
  i) *If $\phi$ satisfies, for a certain $\tilde{m}$,*

$$\phi_x(\tilde{m}, 1) = 0, \quad \nabla_{x,y}(\phi_x)(\tilde{m}, 1) = 0,$$

*it also satisfies* (H).

ii) *There exist a function* $\Lambda(\phi,\psi)$ *and a point* $\tilde{m}(\phi,\psi)$ *such that*
$$\Lambda(\phi^{(0)},0) = 0, \ \tilde{m}(\phi^{(0)},0) = \tilde{m}_0$$
*and the function* $\phi + \psi + \partial_t\phi\Lambda(\phi,\psi)$ *satisfies* (H) *in* $D_0$ *for the point* $\tilde{M} = (\tilde{m}(\phi,\psi),1)$.

iii) *If $\phi$ already satisfies* (H) *for a point* $\tilde{M} = (\tilde{m},1)$ *close to* $\tilde{M}_0$, *then*
$$\Lambda(\phi,0) = 0, \ \tilde{m} = \tilde{m}(\phi,0).$$

*Proof.* Point i) is clear from the Taylor expansion. Let now
$$G : (\phi,\psi,\tilde{m},\lambda) \mapsto (\partial_x\Phi(\tilde{m},1), \ \partial_x^2\Phi(\tilde{m},1), \ \partial_{xy}^2\Phi(\tilde{m},1)) \equiv (G_1, G_2, G_3)$$
with $\Phi = \phi + \psi + \lambda\partial_t\phi$. The function $G$ is of class $C^1$ from $C^3 \times C^3 \times R^3$ to $R^3$. By construction of $\phi^{(0)}$,
$$G(\phi^{(0)}, 0, \tilde{m}_0, 0) = 0.$$

On the other hand,
$$\partial_\lambda G_1(\phi^{(0)}, 0, \tilde{m}_0, 0) = \partial_{xt}^2 \phi^{(0)}(\tilde{m}_0, 1),$$
$$\partial_m G_1(\phi^{(0)}, 0, \tilde{m}_0, 0) = (0,0),$$
$$\partial_m(G_2, G_3)(\phi^{(0)}, 0, \tilde{m}_0, 0) = \nabla^2(\phi_x^{(0)})(\tilde{m}_0, 1) >> 0.$$

Hence the implicit function theorem yields $\lambda = \Lambda(\phi,\psi)$ and $\tilde{m} = \tilde{m}(\phi,\psi)$ with the desired properties. Thanks to i), $\phi + \psi + \partial_t\phi\Lambda(\phi,\psi)$ satisfies (H).

Finally, under the assumptions of iii), $G(\phi, 0, \tilde{m}, 0) = 0$; hence
$$\Lambda(0,0) = 0, \ \tilde{m}(\phi,0) = \tilde{m}.$$

3.4. *Back to the linearized operator.* We go back to Section 3.2 and explain now how we can solve the linearized operator and get $\phi + \dot{\phi}$ to satisfy (H). Assume that $\phi$ already satisfies (H) for $\tilde{m}$ close to $\tilde{m}_0$ and $|\phi - \phi^{(0)}|_{C^4(D_0)}$ small. We will have

(3.3.7) $$\phi + \dot{\phi} = \phi + \dot{\Psi} + \dot{\lambda}\partial_t\phi\frac{\partial_\lambda T}{\partial_t T}.$$

We now take

(3.3.8) $$\dot{\lambda} = (1+\lambda)\Lambda(\phi,\dot{\Psi}) = (1+\lambda)(\Lambda(\phi,\dot{\Psi}) - \Lambda(\phi,0)).$$

Because $\dot{f}$ is small, $\dot{\Psi}$ and $\dot{\lambda}$ are also small: the right-hand side $\hat{\phi}$ of (3.3.7) is then close to $\phi^{(0)}$ and satisfies at $\hat{m} = \tilde{m}(\phi,\dot{\Psi})$
$$\hat{\phi}_x(\hat{m},1) = 0, \ \nabla(\hat{\phi}_x)(\hat{m},1) = 0.$$

According to point i) of the lemma, $\hat{\phi}$ satisfies (H).

**4. An iteration scheme for the problem.** To solve the problem $\mathcal{L}(\lambda, \phi) = 0$ in $D_0$, we will use a Nash-Moser scheme. We refer to [6] for notation and details, and specify here only the nonstandard points.

**4.1. Spaces and smoothing operators.** We will work with the usual Sobolev spaces $H^s(D_0)$. In the process of solving, we note that our starting function $\phi^{(0)}$ satisfies the good boundary conditions, so that all the modifications $\dot\phi$ we will have to consider will be "flat" on $\{t = 0\}$ and $\{x = 0\}$. Hence the smoothing operators used have to respect this "flatness". To achieve this, we take a smooth function $\psi$ supported in $\{t \geq 0, x \leq 0\}$ whose Fourier transform vanishes at the origin of order $k$. Setting

$$S_\theta = \psi_{\theta^{-1}} *, \quad \psi_\varepsilon = \varepsilon^{-3} \psi(x\varepsilon^{-1}, y\varepsilon^{-1}, t\varepsilon^{-1}),$$

we see that the operators $S_\theta$ satisfy the usual properties:
  i) $|S_\theta u|_s \leq C|u|_{s'}$, $s \leq s'$,
  ii) $|S_\theta u|_s \leq C\theta^{s-s'}|u|_{s'}$, $s \geq s'$,
  iii) $|u - S_\theta u|_s \leq C\theta^{s-s'}|u|_{s'}$, $s \leq s'$,
  iv) $|\frac{d}{d\theta} S_\theta u|_s \leq C\theta^{s-s'-1}|u|_{s'}$, but only for $0 \leq s, s' \leq s_k$, $s_k$ going to infinity when $k$ goes to infinity. Here, $|.|_s$ denotes the $H^s$-norm in $D_0$.

**4.2. Smoothing of $\phi$.** In a Nash-Moser procedure, instead of solving

$$\partial_\phi \mathcal{L}(\lambda, \phi)\dot\Psi = \dot f$$

at a given step, we solve

$$\partial_\phi \mathcal{L}(\lambda, \tilde\phi)\dot\Psi = \dot f$$

for an appropriate smoothing $\tilde\phi$ of $\phi$. Wishing to have $\tilde\phi$ satisfy condition (H), we use the following twin of the fundamental Lemma III.3.2.

LEMMA III.3.2′. *Assume $|\phi - \phi^{(0)}|_{C^4(D_0)}$ small enough. Then there exist functions $\tilde\Lambda(\phi)$ and $\tilde m(\phi)$ such that*

$$\tilde\Lambda(\phi^{(0)}) = 0, \quad \tilde m(\phi^{(0)}) = \tilde m_0$$

*and $\phi + x\tilde\Lambda(\phi)$ satisfies* (H) *for the point $\tilde M = (\tilde m, 1)$. Moreover, if $\phi$ satisfies* (H) *for $\tilde m$ close enough to $\tilde m_0$,*

$$\tilde\Lambda(\phi) = 0, \quad \tilde m(\phi) = \tilde m.$$

**4.3. Approximation scheme.** Assuming that we can solve the equation $\partial_\phi \mathcal{L}\dot\Psi = \dot f$ in flat functions, with a tame estimate (see Propositions IV.3.2 and IV.4 for a precise statement), we set

$$\theta_n = (\theta_0^{\eta^{-1}} + n)^\eta, \quad S_n = S_{\theta_n},$$
$$\phi^{(n+1)} = \phi^{(n)} + \Delta\phi^{(n)}, \quad \lambda^{(n+1)} = \lambda^{(n)} + \Delta\lambda^{(n)}.$$

Here, $\Delta$ means "modification" and has nothing to do with the Laplacian! The parameters $\theta_0$ and $\eta^{-1}$ are chosen big enough. For the special smoothing of $\phi^{(n)}$ discussed in 4.2, we set

$$\tilde{S}_n \phi^{(n)} = S_n \phi^{(n)} + x \tilde{\Lambda}(S_n \phi^{(n)}) = S_n \phi^{(n)} + x(\tilde{\Lambda}(S_n \phi^{(n)}) - \tilde{\Lambda}(\phi^{(n)})).$$

Knowing $\lambda^{(n)}, \phi^{(n)}$, we solve in flat functions

$$\partial_\phi \mathcal{L}(\lambda^{(n)}, \tilde{S}_n \phi^{(n)}) \dot{\Psi}_n = \gamma_n$$

for $\gamma_n$ to be determined. Then we take

$$\Delta \lambda^{(n)} = (1 + \lambda^{(n)})[\Lambda(\phi^{(n)}, \dot{\Psi}_n) - \Lambda(\phi^{(n)}, 0)],$$

$$\Delta \phi^{(n)} = \dot{\Psi}_n + \Delta \lambda^{(n)} \partial_t \phi^{(n)} \frac{\partial_\lambda T}{\partial_t T}(\lambda^{(n)}, t).$$

We now determine the $\gamma_n$. First we define the three errors of the solving process:

i) The Taylor error is

$$e'_k = \mathcal{L}(\lambda^{(k+1)}, \phi^{(k+1)}) - \mathcal{L}(\lambda^{(k)}, \phi^{(k)}) - \mathcal{L}'(\lambda^{(k)}, \phi^{(k)})(\Delta \lambda^{(k)}, \Delta \phi^{(k)}).$$

ii) The substitution error is

$$e''_k = [\partial_\phi \mathcal{L}(\lambda^{(k)}, \phi^{(k)}) - \partial_\phi \mathcal{L}(\lambda^{(k)}, \tilde{S}_k \phi^{(k)})] \dot{\Psi}_k.$$

iii) The result error is

$$e'''_k = \Delta \lambda^{(k)} \frac{\partial_\lambda T}{\partial_t T}(t, \lambda^{(k)}) \partial_t(\mathcal{L}(\lambda^{(k)}, \phi^{(k)})).$$

Then we see that

$$\mathcal{L}(\lambda^{(n+1)}, \phi^{(n+1)}) - \mathcal{L}(\lambda^{(n)}, \phi^{(n)}) = e_n + \gamma_n,$$

where $e_n = e'_n + e''_n + e'''_n$ is the total error. Finally, we denote by

$$E_n = \Sigma_{0 \le k \le n-1} e_k$$

the accumulated error, and

$$\mathcal{L}(\lambda^{(n+1)}, \phi^{(n+1)}) = f^{(0)} + \Sigma_{0 \le k \le n} \gamma_k + E_{n+1}.$$

It is natural at this stage to determine the $\gamma_n$ by

$$\Sigma_{0 \le k \le n} \gamma_k + S_n E_n = -S_n f^{(0)},$$

which leads to

$$\mathcal{L}(\lambda^{(n+1)}, \phi^{(n+1)}) = f^{(0)} - S_n f^{(0)} + e_n + (E_n - S_n E_n).$$

### IV. Existence and tame estimates for the linearized problem

1. *Structure of the linearized operator.* In part III, we showed how solving $\mathcal{L}(\lambda, \phi) = 0$ can be reduced to solving the linearized equation $\partial_\phi \mathcal{L}$. We display now the structure of this operator.

PROPOSITION IV.1. *The linearized operator has the form*

(4.1.1) $$c_0 \partial_\phi \mathcal{L}(\lambda, \phi) = \hat{Z}\bar{S}\hat{Z} + \varepsilon^2(\bar{S}\phi)\hat{N}\hat{Z} + \varepsilon^2 B \partial_y^2 + b_0 \bar{S}\hat{Z} + \varepsilon^2 \ell.$$

*Here*

i) $\hat{Z}$ *is a field of the form*
$$\hat{Z} = \partial_t + \varepsilon^2 z_0 \partial_y, \quad z_0 = z_0(x, y, t, \lambda, \phi, \phi_y, \phi_t),$$

ii) $\bar{S}$ *is the field (independent of $\phi$)*
$$\bar{S} = \partial_x + \varepsilon^2 s_0 \partial_t, \quad s_0 = s_0(x, y, t, \lambda),$$

iii) $\hat{N}$ *is the second order operator*
$$\hat{N} = \hat{N}_1 \hat{Z}^2 + 2\varepsilon^2 \hat{N}_2 \hat{Z}\partial_y + \hat{N}_3 \partial_y^2, \quad \hat{N}_i = \hat{N}_i(x, y, t, \lambda, \phi, \phi_y, \phi_t),$$
$$\hat{N}_1 = -\frac{1}{4\partial_t T(\tau_0 + T(t, \lambda))} + O(\varepsilon^2), \quad \hat{N}_3 = \frac{\partial_t T}{(\tau_0 + T(t, \lambda))^3} + O(\varepsilon^2),$$

iv) $B = -(\hat{Z}\bar{S}\phi)\hat{N}_3$,

v) $c_0 = c_0(x, y, t, \lambda, \phi, \phi_y, \phi_t) = (q\partial_t T)^{-2} + O(\varepsilon^2),$
$b_0 = b_0(x, y, t, \lambda, \nabla\phi, \nabla^2\phi),$

vi) $\ell$ *is a second order operator which can be written as a linear combination of*
$$\text{id}, \bar{S}, \hat{Z}, \partial_y, \bar{S}\hat{Z}, \hat{Z}^2, \hat{Z}\partial_y, \partial_y^2$$
*with coefficients depending on the derivatives of $\phi$ up to order 3. Moreover, $\ell$ does not contain $\partial_y^2$ for $\varepsilon = 0$.*

*Proof.* a. The linearized operator $\partial_\phi \mathcal{L}$ is obtained as follows: first, we linearize the $\phi$ equation resulting from substituting (3.2.4) into (3.1.16). Then we perform the changes of variables (3.2.3) and (3.3.1).

b. With the notation of Lemma III.1, let us compute $\dot{\mathcal{T}}_0$. We find
$$\dot{\mathcal{T}}_0 = q\dot{v} + Z\dot{\phi} + \varepsilon^2 \gamma_1 \dot{\phi}, \quad \gamma_1 = \gamma_1(\omega, \tau, \phi, \phi_\omega, \phi_\tau).$$
Hence $\dot{V} = -\frac{1}{q}(Z\dot{\phi} + \varepsilon^2 \gamma_1 \dot{\phi})$. On the other hand, linearizing $\mathcal{T}_2 + \phi_s \mathcal{T}_3 = 0$ gives
$$\dot{Z}V_s + Z\dot{V}_s + \varepsilon^2 \dot{\phi}_s NV + \varepsilon^2 \phi_s \dot{N}V + \varepsilon^2 \phi_s N\dot{V} - \varepsilon^2 \dot{V}_s N\phi$$
$$- \varepsilon^2 V_s \dot{N}\phi - \varepsilon^2 V_s N\dot{\phi} + \varepsilon^2 \dot{V}_s h_2 + \varepsilon^2 V_s \dot{h}_2 + \varepsilon^2 \dot{\phi}_s h_3 + \varepsilon^2 \phi_s \dot{h}_3.$$

We see that $\dot{Z}$ and $\dot{N}$ yield only first order derivatives of $\dot{\phi}$, multiplied by $\varepsilon^2$. The same is true of $\dot{h}_2$ and $\dot{h}_3$, except for the terms $\dot{V}_\omega$ and $\dot{V}_\tau$ coming from the corresponding terms $v_\omega$ and $v_\tau$ in $h_2$ and $h_3$. It follows that the linearized equation on $\dot{\phi}$ has the form

$$\{Z\partial_s + \varepsilon^2\phi_s N - \varepsilon^2(N\phi)\partial_s + \varepsilon^2 h_2\partial_s + \varepsilon^2 h_4\partial_\omega$$
$$+ \varepsilon^2 h_5\partial_\tau\}\dot{V} - \varepsilon^2 V_s N\dot{\phi} + \varepsilon^2 \times \nabla\dot{\phi},$$

the last term denoting just first order derivatives of $\dot{\phi}$.

c. The composition of the two changes of variables, denoted by bars, operates the following transformation of operators (to avoid introducing unnecessary notation, we denote by $*$ known functions):

$$\bar{\partial}_s = \partial_x + \varepsilon^2 s_0(x, y, t, \lambda)\partial_t \equiv \bar{S},$$
$$\bar{\partial}_\omega = \partial_y + \varepsilon^2 * (x, y, t, \lambda)\partial_t,$$
$$\bar{\partial}_\tau = (1 + \varepsilon^2 * (x, y, t, \lambda))(\partial_t T)^{-1}\partial_t,$$
$$\bar{Z} = (-1 + \varepsilon^2 * (x, y, t, \lambda, \phi, \phi_y, \phi_t))\hat{Z},$$
$$\bar{N} = \bar{N}^{(1)}\partial_t^2 + 2\varepsilon^2 \bar{N}^{(2)}\partial_{yt}^2 + \bar{N}^{(3)}\partial_y^2 + \text{lower order terms},$$
$$\bar{N}^{(1)} = \frac{1}{4(\partial_t T)^2(\tau_0 + T(t, \lambda))} + \varepsilon^2 * (x, y, t, \lambda, \phi, \phi_y, \phi_t),$$
$$\bar{N}^{(3)} = -\frac{1}{(\tau_0 + T(t, \lambda))^3} + \varepsilon^2 * (x, y, t, \lambda, \phi, \phi_y, \phi_t).$$

Finally, we replace the operators by their transforms into the linearized equation and set $\hat{N} = -(\partial_t T)\bar{N}$. To obtain the value of $B$, we observe, keeping only the $\varepsilon^2$ terms in the coefficient of $\partial_y^2$, that

$$B = \bar{q}(\partial_t T)^2(\bar{S}\bar{V})\bar{N}^{(3)} + O(\varepsilon^2),$$

from which iv) follows. $\square$

2. *Energy inequality for the linearized operator.* In the following, in a (desperate) attempt to simplify the notation, we will write abusively $Z$ for $\hat{Z}$, $S$ for $\bar{S}$, $N$ for $\hat{N}$, and replace $\varepsilon^2$ everywhere by $\varepsilon$.

We assume a given smooth function $\phi$ in $D_0$, close to $\phi^{(0)}$, satisfying (H) for a point $\tilde{M} = (\tilde{m}, 1)$. We then set

$$P \equiv ZSZ + \varepsilon(S\phi)NZ + \varepsilon B\partial_y^2 + \varepsilon\ell + b_0 SZ,$$
$$\tilde{P} = ZSZ + \varepsilon(S\phi)NZ.$$

Recall that we have arranged for $\{t = 0\}$ to be characteristic, that is

(4.2.1) $\qquad (t = 0) \Rightarrow s_0 + (S\phi)N_1 = 0.$

The connection between $P$ and $\tilde{P}$ is explained in the following straightforward lemma.

LEMMA IV.2. *With the above notation,*

$$ZPu = \tilde{P}Zu + \varepsilon\tilde{\ell}Zu + \varepsilon Z\ell u + Z(b_0 SZu) + \varepsilon[ZB - 2\varepsilon B(\partial_y z_0)]\partial_y^2 u - \varepsilon^2 B(\partial_y^2 z_0)\partial_y u.$$

Here, $\tilde{\ell}$ *is a second order operator which can be written as a linear combination of* $Z^2$, $Z\partial_y$, $\partial_y^2$, $\partial_y$, *the coefficient of* $\partial_y^2$ *being a multiple of* $(S\phi)$.

The point of this lemma is that $\tilde{P}$ does not contain the delicate term $\varepsilon B\partial_y^2$, as $P$ does.

2.1. *Energy inequality for* $\tilde{P}$. The following computation is very close to that of [4, Section III.2]. Unfortunately, $N$ is different here, causing some new problems. Thus we think it better to give the whole computation again, which is, after all, the heart of the proof. We set

$$A = S\phi, \quad \delta = 1 - t, \quad g = \exp h(x - t),$$
$$p = \delta^{\frac{\mu}{2}} \exp \frac{h}{2}(x - t), \quad |\cdot|_0 = |\cdot|_{L^2(D_0)}.$$

PROPOSITION IV.2.1. *Fix* $\mu > 1$. *Then there exist* $C > 0$, $\varepsilon_0 > 0$, $\eta_0 > 0$ *and* $h_0$ *such that, for all* $\phi$ *satisfying* (H), (4.2.1) *and* $|\phi - \phi^{(0)}|_{C^4(D_0)} \leq \eta_0$, $0 \leq \varepsilon \leq \varepsilon_0$, $h \geq h_0$ *and smooth* $u$ *with*

$$u(x, y, 0) = \partial_t u(x, y, 0) = 0, \quad u(M, y, t) = 0,$$

*the following inequality exists*:
(4.2.2)
$$h|pSZu|_0^2 + h|pZ^2u|_0^2 + \varepsilon h|p\partial_y Zu|_0^2 + \varepsilon^2 \int \delta^{\mu-1} g(S\phi)(1+\delta h)|\partial_y^2 u|^2 \leq C|p\tilde{P}u|_0^2.$$

*Proof.* a. With a still unknown multiplier

$$Mu = aSZu + \varepsilon c\partial_y^2 u + dZ^2 u,$$

we find by the usual integrations by parts

$$\int \tilde{P}uMu\,dxdydt$$
$$= \int K_0(SZu)^2 + \int \varepsilon^2 K_1(\partial_y^2 u)^2 + \int \varepsilon K_2(\partial_y Zu)^2 + \int K_3(Z^2u)^2$$
$$+ \varepsilon \int (\partial_y Zu)(SZu)[\partial_y(c - aAN_3) - \varepsilon Z(aAN_2) - 2\varepsilon^2 aAN_2\partial_y z_0]$$
$$+ \varepsilon \int (\partial_y^2 u)(SZu)[-Zc + \varepsilon c\partial_y z_0]$$
$$+ \varepsilon \int (Z^2u)(\partial_y Zu)[-\partial_y(dAN_3) + \varepsilon S(aAN_2) + \varepsilon\partial_y(cAN_1)$$
$$+ \varepsilon aA(N_1 Sz_0 - N_3\partial_y s_0 - \varepsilon^2 z_0\partial_y s_0) + \varepsilon c\partial_y s_0 - 2\varepsilon^2 dAN_2\partial_y z_0]$$

$$+ \varepsilon^2 \int (Z^2 u)(\partial_y^2 u)[-Z(cAN_1) + \varepsilon cAN_1 \partial_y z_0]$$
$$+ \varepsilon^3 \int (\partial_y^2 u)(\partial_y Zu)[-Z(cAN_2) + \varepsilon cAN_2 \partial_y z_0]$$
$$+ \varepsilon \int (Z^2 u)(SZu)[-Z(aAN_1) - \varepsilon \partial_y(aAN_2) - \varepsilon aAN_1 \partial_y z_0]$$
$$+ \varepsilon^2 \int (\partial_y u)(\partial_y^2 z_0)[cSZu + cAN_1 Z^2 u + cAN_3 \partial_y^2 u + \varepsilon cAN_2 \partial_y Zu]$$
$$+ I_1 - I_0 + J_1.$$

The coefficients $K_i$ of the quadratic terms are

$$2K_0 = -Za - \varepsilon a \partial_y z_0,$$
$$2K_1 = -Z(cAN_3) + 3\varepsilon cAN_3 \partial_y z_0,$$
$$2K_2 = S(aAN_3 - c) + Z(dAN_3) - \varepsilon Z(cAN_1) + \varepsilon^2 \partial_y(cAN_2) - \varepsilon^2 dAN_3$$
$$\quad + \varepsilon^2 aAN_3(Zs_0 + \varepsilon z_0 \partial_y s_0) - c\varepsilon^2(Zs_0 + \varepsilon z_0 \partial_y s_0 - \varepsilon AN_1 \partial_y z_0)$$
$$\quad + 2\varepsilon^3 aAN_2(Sz_0 + \varepsilon z_0 Zs_0),$$
$$2K_3 = -Sd + \varepsilon S(aAN_1) - \varepsilon Z(dAN_1)$$
$$\quad - \varepsilon d(\partial_t s_0 + \varepsilon AN_1 \partial_y z_0) - 2\varepsilon^2 \partial_y(dAN_2)$$
$$\quad + \varepsilon^2 aA(N_1 \partial_t s_0 - 2N_1 Z s_0 - 2\varepsilon N_2 \partial_y s_0).$$

The terms $I_0, I_1, J_1$ are boundary terms. We have

$$2I_1 = \int_{\{t=1\}} [a(SZu)^2 + \varepsilon dAN_1(Z^2 u)^2$$
$$\quad - \varepsilon(dAN_3 - \varepsilon cAN_1)(\partial_y Zu)^2 + \varepsilon^2 cAN_3(\partial_y^2 u)^2$$
$$\quad + 2\varepsilon(SZu)(aAN_1(Z^2 u) + \varepsilon aAN_2 \partial_y Zu + c\partial_y^2 u)$$
$$\quad + 2\varepsilon^2 cAN_1(Z^2 u)(\partial_y^2 u) + 2\varepsilon^3 cAN_2(\partial_y Zu)(\partial_y^2 u)].$$

The $I_0$ term on $\{t = 0\}$ vanishes, thanks to (4.2.1) and to the fact that the first two traces of $u$ are zero. Finally,

$$2J_1 = \int_{\{x=-A_0\}} (-d + \varepsilon aAN_1)(Z^2 u)^2 - \varepsilon(c - aAN_3)(\partial_y Zu)^2$$
$$\quad + 2\varepsilon^2 aAN_2(\partial_y Zu)(Z^2 u).$$

b. We choose now

$$a = A^{-1} \delta^\mu g, \quad c = c' \delta^\mu g, \quad d = -d' \delta^\mu g$$

where $h, c', d'$ are positive constants to be chosen later. We analyze first the boundary terms.

We find that $2I_1$ is the integral of

$$a(SZu + \varepsilon AN_1 Z^2 u + \varepsilon^2 AN_2 \partial_y Zu + c'A\partial_y^2 u)^2$$

$$- \varepsilon(Z^2u)^2\delta^\mu g A N_1(d' + \varepsilon N_1) + \varepsilon(\partial_y Zu)^2 \delta^\mu g A(d'N_3 + \varepsilon c'N_1 - \varepsilon^3 N_2^2)$$
$$+ \varepsilon^2(\partial_y^2 u)^2 \delta^\mu g c' A(N_3 - c') - 2\varepsilon^3 A\delta^\mu g N_1 N_2(Z^2u)(\partial_y Zu).$$

From Proposition IV.1, it is clear that $I_1 \geq 0$ for $c' < N_3$ and $\varepsilon$ small enough.

The term $2J_1$ is the integral of

$$(\delta^\mu g)[(d' + \varepsilon N_1)(Z^2u)^2 + \varepsilon(N_3 - c')(\partial_y Zu)^2 + 2\varepsilon^2 N_2(\partial_y Zu)(Z^2u)].$$

Clearly, $J_1 \geq 0$ for $c' < N_3$ and $\varepsilon$ small enough.

c. We analyze now the signs of the quadratic terms.

We find

$$2K_0 = \frac{\delta^{\mu-1}g}{A^2}[\delta ZA + \mu A + \delta A(h - \varepsilon \partial_y z_0)].$$

In a small neighbourhood $\omega$ of $M$, we have $A = \phi_x \geq -\delta V$, with

$$V \leq 0, \quad V(M) = \phi_{xt}(M) < 0$$

according to (H). It follows that, with $\mu' = \frac{\mu+1}{2}$ and a possibly smaller $\omega$,

$$\delta ZA + \mu'A \geq \delta(\phi_{xt} + \varepsilon z_0 \phi_{xy} - \mu'V) \geq 0;$$

hence

$$K_0 \geq C\frac{\delta^{\mu-1}g}{A}(1 + \delta h).$$

Outside $\omega$, we have for $h$ big enough

$$\delta(ZA + A(h - \varepsilon \partial_y z_0)) \geq 1/2\delta Ah;$$

so finally

$$K_0 \geq k_0 \frac{\delta^{\mu-1}g}{A}(1 + \delta h), \quad k_0 > 0.$$

Now

$$2K_1 = c'\delta^{\mu-1}g A N_3(\mu + 3\varepsilon\delta\partial_y z_0) + c'\delta^\mu g[(hN_3 - ZN_3)A - (ZA)N_3].$$

Close to $M$, $ZA < 0$, hence

$$(hN_3 - ZN_3)A - (ZA)N_3 \geq ChA.$$

This also holds away from $M$ for big $h$. Finally,

$$K_1 \geq k_1 c'\delta^{\mu-1}g A(1 + \delta h).$$

The first two terms in $2K_2$ are

$$\delta^\mu g(SN_3) + g(N_3 - c')(h\delta^\mu(1 - \varepsilon s_0) - \varepsilon\mu s_0 \delta^{\mu-1})$$
$$+ \mu\delta^{\mu-1}d'gAN_3 + d'\delta^\mu g[(hN_3 - ZN_3)A - (ZA)N_3].$$

For $c' < N_3$, they are bigger than

$$2k_2 g(\delta^\mu h + \delta^{\mu-1}A).$$

All the other terms in $K_2$ are bounded by
$$C\varepsilon\delta^{\mu-1}gA + C\varepsilon h\delta^\mu g.$$

Hence
$$K_2 \geq k_2\delta^{\mu-1}g(A + \delta h).$$

Finally,
$$\begin{aligned}2K_3 &= g(d' + \varepsilon N_1)[h\delta^\mu(1 - \varepsilon s_0) - \varepsilon\mu s_0\delta^{\mu-1}] \\ &\quad + \varepsilon\delta^\mu g(SN_1) - \varepsilon N_1 d'A\delta^{\mu-1}g(\mu + h\delta) \\ &\quad + O(\varepsilon\delta^\mu g) \geq k_3 g(h\delta^\mu + \varepsilon\delta^{\mu-1}A).\end{aligned}$$

d. Consider now the product term in $(\partial_y^2 u)(SZu)$. To ensure positivity, it is enough to check the positivity of the quadratic form
$$-\theta Za(SZu)^2 - \varepsilon^2\theta Z(cAN_3)(\partial_y^2 u)^2 - 2\varepsilon(Zc)(\partial_y^2 u)(SZu)$$
for some $\theta < 1$. The discriminant $\Delta$ satisfies
$$-\varepsilon^{-2}\Delta = \theta^2(Za)(Z(cAN_3)) - (Zc)^2.$$

After rearrangement of terms we get
$$\begin{aligned}-\frac{c'}{c^2}\delta^2 A^2\varepsilon^{-2}\Delta &= (\theta^2 N_3 - c')A^2(\mu + h\delta)^2 - \theta^2 N_3\delta^2(ZA)^2 \\ &\quad - \theta^2 A^2(ZN_3)\delta(\mu + h\delta) - \theta^2\delta^2 A(ZA)(ZN_3).\end{aligned}$$

Away from $M$, $A$ is uniformly positive and the right-hand side is positive for big $h$ if $c' < \theta^2 N_3$. Close to $M$,
$$A \geq -\delta V, \ V(M) = ZA - \varepsilon z_0\partial_y A,$$
and the right-hand side is bigger than
$$[(\theta^2 N_3 - c')\mu^2 + O(\delta)]\delta^2 V^2 - \theta^2 N_3\delta^2(ZA)^2.$$

For $c'$ small enough and $\theta$ close to 1, this is bigger than $C\delta^2$.

e. It is easy (but lengthy) to check that all other product terms can be absorbed in the quadratic terms for small $\varepsilon$.

f. Finally, using Cauchy-Schwarz inequality, we obtain (4.2.2). □

2.2. *Energy inequality for $P$.* We deduce from (4.2.2) and Lemma IV.2 an energy inequality for $P$.

PROPOSITION IV.2.2. *Fix $\mu > 1$. Then there exist $C > 0$, $\epsilon_0 > 0$, $\eta_0 > 0$ and $h_0$ such that, for all $\phi$ satisfying* (H), (4.2.1) *and $|\phi - \phi^{(0)}|_{C^4(D_0)} \leq \eta_0$, $0 \leq \varepsilon \leq \varepsilon_0$, $h \geq h_0$ and smooth $u$ with*
$$u(x,y,0) = \partial_t u(x,y,0) = 0, \ u(M,y,t) = 0, \ Pu(x,y,0) = 0,$$

*there exists the energy inequality*

(4.2.4) $$h|pSZ^2u|_0^2 + h|pZ^3u|_0^2 + \varepsilon h|pZ^2\partial_y u|_0^2 + \varepsilon^2|pZ\partial_y^2 u|_0^2$$
$$+ \varepsilon^2 h \int \delta^\mu g(S\phi)|\partial_y^2 Zu|^2 \leq C|pZPu|_0^2.$$

*Proof.* a. First, let us compute in $Pu$ the coefficient of $\partial_t^2 u$ on $\{t = 0\}$. We get
$$[\partial_x + \varepsilon^2(2s_0 z_0 + 2AN_2 + 3Az_0)\partial_y + \varepsilon^2 z_0 \partial_y s_0](\partial_t^2 u).$$

Hence the assumptions on $u$ imply $\partial_t^2 u(x,y,0) = 0$.

b. We can then apply inequality (4.2.2) to $Zu$, getting
$$h|pSZ^2u|_0^2 + h|pZ^3u|_0^2 + \varepsilon h|p\partial_y Z^2 u|_0^2$$
$$+ \varepsilon^2 \int \delta^{\mu-1} g(S\phi)(1+\delta h)|\partial_y^2 Zu|^2 \leq C|p\tilde{P}Zu|_0^2.$$

We now have to bound the terms of $\tilde{P}Zu - ZPu$ (given in Lemma IV.2) by the left-hand side. Note first that, by standard lemmas (see for instance [2, Lemma 2.1]), the terms of the left-hand side also give control of
$$\varepsilon h|pZ^2\partial_y u|_0^2 + \varepsilon h|pZ\partial_y Zu|_0^2 + h|pZSZu|_0^2.$$

Thus the terms from $\varepsilon(\tilde{\ell} + \ell)Zu$ are easily dominated for small $\varepsilon$ and big $h$.

Since $\delta \leq C\phi_x$, we have control of $\varepsilon^2|p\partial_y^2 Zu|_0^2$, implying control of both $\varepsilon^2|pZ\partial_y^2 u|_0^2$ and $\varepsilon^2 h|p\partial_y^2 u|_0^2$. This makes it possible to absorb the terms $\varepsilon[Z,\ell]$ and also the remaining terms in $\partial_y^2 u$, $\partial_y u$, $SZu$ and $ZSZu$ for small $\varepsilon$ and big $h$. □

3. *Higher order inequalities.*

3.1. *The spaces $\tilde{H}^s$.* In Section 2, we used the fields $Z, S, \partial_y$ systematically instead of the standard set $\partial_t, \partial_x, \partial_y$. The reason for this is that if we develop the expression of $P$, the energy inequality becomes much less transparent (to say the least); in fact, one should observe that (4.2.2) by itself does not give separate control of $\partial_{xy}^2 u$, for instance. Thus it is also appropriate to commute $P$ with $Z, S, \partial_y$ and their products. We define first $T^l$ as any product of $l$ fields among $Z, S, \partial_y$, and for integer $s$,
$$\tilde{H}^s = \{u \in L^2(D_0), l \leq s \Rightarrow T^l u \in L^2(D_0)\}.$$

We denote the natural norm by
$$|u|_{\tilde{s}}^2 = \Sigma |T^l u|_0^2.$$

PROPOSITION IV.3.1. *The spaces $H^s$ and $\tilde{H}^s$ are the same. Moreover, for $\phi$ bounded in $C^1$,*

i) $$|u|_{\tilde{s}}^2 \leq C(|u|_s^2 + |u|_{L^\infty}^2(1 + |\phi|_{s+1}^2)),$$
ii) $$|u|_s^2 \leq C(|u|_{\tilde{s}}^2 + |u|_{L^\infty}^2(1 + |\phi|_{s+1}^2)).$$

*Proof.* a. The equality of the spaces is obvious for smooth $\phi$.

b. We will repeatedly use the following classical lemma (see for instance [6]).

LEMMA. *For $|\alpha_1| + \ldots + |\alpha_p| = s$,*
$$|\partial^{\alpha_1} u_1 \ldots \partial^{\alpha_p} u_p|_0 \leq C \Sigma_{1 \leq j \leq p} |u_1|_{L^\infty} \ldots |u_j|_s \ldots |u_p|_{L^\infty}.$$

Denoting by $\partial$ usual derivatives, we have (skipping everywhere irrelevant numerical coefficients)
$$T^l u = \Sigma \partial^{q_1} a \ldots \partial^{q_k} a \partial^r u,$$
where $a = \varepsilon s_0$ or $\varepsilon z_0$, and
$$r \geq 1, k \leq l, |q_1| + \ldots + |q_k| + |r| = l.$$
Hence
$$|T^l u|_0 \leq C(|u|_l + |u|_{L^\infty} |a|_l),$$
which gives i).

c. Conversely, we have
$$\partial_t = Z - \varepsilon z_0 \partial_y, \quad \partial_x = S - \varepsilon s_0 Z + \varepsilon^2 s_0 z_0 \partial_y.$$
For $p + q + m = l$, we can write
$$\partial_t^p \partial_x^q \partial_y^m = \Sigma T^{q_1}(a) \ldots T^{q_k}(a) T^r u,$$
with
$$|q_1| + \ldots + |q_k| + |r| = p + q + m, \quad k \leq p + q + m$$
and $a = \varepsilon z_0, \varepsilon s_0$ or $\varepsilon^2 s_0 z_0$. Using the identity of b, except for the terms $T^l u$, we get
$$\partial_t^p \partial_x^q \partial_y^m = \Sigma a \ldots a T^l u + \Sigma \partial^{p_1} a \ldots \partial^{p_j} a \partial^r u,$$
with $|p_1| + \ldots + |p_j| + |r| = l, r \leq l - 1$ in the last sum. Taking $L^2$ norms and using the interpolation lemma, one obtains
$$|\partial_t^p \partial_x^q \partial_y^m u|_0 \leq C|T^l u|_0 + C|u|_{l-1} + C|u|_{L^\infty} |a|_l.$$
Hence, by induction on $s$, we have ii). □

3.2. *Structure of the commutators.* The two following straightforward lemmas describe the structure of the commutators of $P$ with a product $K = T^l$.

LEMMA IV.3.1. *Let us denote here by $c$ various coefficients which are smooth functions of*
$$x, y, t, \lambda, \phi, \nabla\phi, \nabla^2\phi,$$
*and by $Q_3$ the principal part of the third order operator $Q$. Then*
 i) $[Z, P]_3 = \varepsilon(cZ^3 + c\partial_y Z^2 + c\partial_y^2 Z)$,
 ii) $[S, P]_3 = \varepsilon(cZ^3 + c\partial_y Z^2 + c\partial_y^2 Z + cZSZ + c\partial_y SZ + c(S\phi)\partial_y^3)$,
 iii) $[\partial_y, P]_3 = \varepsilon(cZ^3 + c\partial_y Z^2 + c\partial_y^2 Z + c(S\phi)\partial_y^3)$.

LEMMA IV.3.2. *When $K = T^l$,*
$$Z[K, P]u = Z\Sigma[T, P]_3 T^{l-1} u + Z\Sigma\partial^{q_1} a \ldots \partial^{q_k} a \partial^q b T^r u$$
$$+ \Sigma\partial^{q_1} a \ldots \partial^{q_k} a \partial^q b T^l u + \Sigma\partial^{q_1} a \ldots \partial^{q_k} a \partial^q b \partial^r u$$
$$\equiv Z\Sigma_1 + Z\Sigma_2 + \Sigma_3 + \Sigma_4.$$

Here, $a = \varepsilon s_0$ or $\varepsilon z_0$, $b$ is a coefficient of $P$, and the conditions on the derivatives are the following:
 i) In $\Sigma_2$, $l \leq r \leq l+1$, $|q| + |q_1| + \ldots + |q_k| + |r| \leq l+3$.
 ii) In $\Sigma_3$, $|q| + |q_1| + \ldots + |q_k| \leq 4$.
 iii) In $\Sigma_4$, $r \leq l-1, |q| + |q_1| + \ldots + |q_k| + |r| \leq l+4$.

*Moreover, all the terms in $Z[K, P]u$ contain a factor $\varepsilon$, except for $Z[K, b_0 SZ]u$ which is of the special form*
$$Z[K, b_0 SZ]u = Z^2 \Sigma \ldots T^l u + ZS\Sigma \ldots T^l u + Z\Sigma \ldots T^l u$$
$$+ \Sigma \ldots T^l u + \Sigma_{r \leq l-1} \ldots \partial^r u.$$

3.3. *Higher order tame estimates.*

PROPOSITION IV.3.1. *Let $\phi$ satisfy the assumptions of Proposition IV.2.1 and moreover $|\phi - \phi^{(0)}|_{C^6(D_0)} \leq \eta_0$. Then there exists $\varepsilon_0$ such that for $0 \leq \varepsilon \leq \varepsilon_0$ and all $s$, there exists $C_s$ such that for all smooth $u$ flat on $\{t = 0\}$ and $\{x = M\}$,*

(4.3.1) $\qquad |u|_s \leq C|\partial_\phi \mathcal{L}u|_{s+1} + C|\partial_\phi \mathcal{L}u|_3 (1 + |\phi|_{s+7}).$

*Proof.* a. Recall that $c_0 \partial_\phi \mathcal{L} = P$. Thus, it is enough to prove (4.3.1) for $P$. Let $Pu = f, K = T^l, l \leq s$. Then
$$KZf = [K, Z]f + Z[K, P]u + ZP(Ku).$$
Applying the energy inequality (4.2.4) to $Ku$ and summing over $l \leq s$, we obtain control of a certain number of terms which we denote by $E_s$ and do not repeat here. We have to show that the terms of $Z[K,P]u$, whose structure is discussed in Lemmas IV.3.1 and IV.3.2, can be bounded in the weighted norm $|p.|_0$ by an arbitrarily small fraction of $E_s$ (in short, "absorbed").

b. Examine first the terms of $Z\Sigma_1$. Taking into account the explicit expressions of the commutators $[T,P]_3$ given in the lemma, we see that there are six different types of terms to control, which we display in three groups:

$$(1) \qquad \varepsilon Z(cZ^3)T^{l-1}u, \ \ \varepsilon Z(c\partial_y Z^2)T^{l-1}u, \ \ \varepsilon Z(cZSZ)T^{l-1}u,$$

$$(2) \qquad \varepsilon Z(c\partial_y^2 Z)T^{l-1}u, \ \ \varepsilon Z(c\partial_y SZ)T^{l-1}u,$$

$$(3) \qquad \varepsilon Z(c(S\phi)\partial_y^3)T^{l-1}u.$$

The terms of group (1) are clearly absorbed for big $h$.

For the group (2), we write

$$\partial_y^2 ZT^{l-1}u = \partial_y Z(\partial_y T^{l-1}u) + *\partial_y(\partial_y T^{l-1}u) + *(\partial_y T^{l-1}u),$$

which shows absorption for big $h$. We proceed analogously by splitting the second term in terms easily absorbed

$$\partial_y SZT^{l-1}u = SZ(\partial_y T^{l-1}u) + *S(\partial_y T^{l-1}u) \\ + *Z^2 T^{l-1}u + *\partial_y ZT^{l-1}u + *\partial_y T^{l-1}u.$$

Finally,

$$Z(c(S\phi)\partial_y^2)(\partial_y T^{l-1}u) = (Zc)(S\phi)\partial_y^2(\partial_y T^{l-1}u) + cZ(S\phi)\partial_y^2(\partial_y T^{l-1}u) \\ + c(S\phi)\partial_y^2 Z(\partial_y T^{l-1}u) + c(S\phi)[Z,\partial_y^2](\partial_y T^{l-1}u).$$

Since $[Z,\partial_y^2] = *\partial_y^2 + *\partial_y$, all the terms can be absorbed for big $h$.

c. We analyze now $Z\Sigma_2$ and $\Sigma_3$.

A product $T^{s+1}$ containing at least one factor $S$ can be written

$$T^{s+1}u = ST^s u + \Sigma \partial a T^s u + \Sigma \partial^{q_1} a \ldots \partial^{q_k} a \partial^r u,$$
$$r \leq s-1, \ \ |r| + |q_1| + \ldots + |q_k| = s+1.$$

We can proceed similarly if $T^{s+1}$ contains at least a factor $Z$. In all cases,

$$\varepsilon h|pT^{s+1}u|_0^2 \leq CE_s + C\varepsilon h|p\Sigma_{r \leq s-1} \ldots \partial^r u|_0^2.$$

A similar analysis gives the same estimate for terms like $\varepsilon h|pZT^{s+1}u|_0^2$. Thus, all terms from $Z\Sigma_2$ or $\Sigma_3$ which have $\varepsilon$ as a factor are absorbed for big $h$, modulo an additional term $\varepsilon h|p\Sigma_{r \leq s-1} \ldots \partial^r u|_0^2$ on the right-hand side of the inequality. The same is true for the terms from $Z[K,b_0 SZ]u$, thanks to their special structure.

d. We now fix $h$ and use the inequalities

$$|pv|_0 \leq C|v|_0, \ \ E_s \geq Ch|u|_{\tilde{s}}.$$

Applying the interpolation lemma for the index $s-1$ ("taking out" five derivatives from the coefficients), we can bound the additional terms and the terms of $\Sigma_4$:

$$|\partial^{q_1} a \ldots \partial^{q_k} a \partial^q b \partial^r u|_0 \leq C|u|_{s-1} + C|u|_{L^\infty}(1 + |\phi|_{s+7}).$$

By induction on $s$, using Proposition IV.3.1, we finally obtain

$$|u|_s^2 \leq C|f|_{s+1}^2 + C|f|_3^2(1 + |\phi|_{s+6}^2) + C|u|_2^2(1 + |\phi|_{s+7}^2).$$

Since it is easy to obtain a low norm estimate

$$|u|_2 \leq C|f|_3,$$

we get (4.3.1). □

4. *Existence of flat solutions.*

PROPOSITION IV.4. *Let $\phi$ satisfy the assumptions of Proposition IV.2.2. Then there exists, for all $f \in C^\infty(D_0)$, flat on $\{t = 0\}$ and $\{x = M\}$, a unique smooth solution of*

$$\partial_\phi \mathcal{L} u = f, u(x, y, 0) = \partial_t u(x, y, 0) = 0, \quad u(M, y, t) = 0,$$

*satisfying the estimates* (4.3.1).

*Proof.* a. Define, as in [4], the smoothed operators $P_\alpha$ by replacing $\partial_y$ by $Y_\alpha = \chi(\alpha D_y)\partial_y$ in $P$ ($\chi$ being compactly supported and one near zero). For fixed $\alpha > 0$, we can solve

$$P_\alpha u = f, \quad u(x, y, 0) = \partial_t u(x, y, 0) = 0, \quad u(M, y, t) = 0$$

in smooth functions. In fact, setting

$$\tilde{S} = S + \varepsilon(S\phi)N_1\partial_t = \partial_x + \varepsilon(s_0 + (S\phi)N_1)\partial_t \equiv \partial_x + \varepsilon\tilde{s}_0\partial_t,$$

we can expand the terms of $P_\alpha$ and write

$$P_\alpha u \equiv \partial_t^2 \tilde{S} u + \partial_t \tilde{S} A_1 u + \partial_t^2 A_2 u + \partial_t A_3 u + \tilde{S} A_4 u + A_5 u.$$

Here, the $A_i$ are zero order operators in $y$ depending smoothly on $(x, t)$. To solve in $D_0$, we proceed as usual by writing the equation in integral form.

More precisely, let $M_0 = (x_0, y_0, t_0)$ be some point in $D_0$ and denote by

$$X(M_0, x), Y(M_0, x), T(M_0, x)$$

the parametrization by $x$ of the segment of the integral curve of $\tilde{S}$ from $M_0$. We set

$$(I_{\tilde{S}} v)(M_0) = \int_M^{x_0} v(X, Y, T) dx,$$

so $I_{\tilde{S}}v = 0$ for $x = M$ and $\tilde{S}I_{\tilde{S}}v = v$. The key point is that this segment of integral curve (for $x_0 \leq x \leq M$) is well defined and contained in $D_0$ because of the following two facts:

i) On $\{t = 0\}$, $\tilde{S} = \partial_x$,

ii) Since $s_0 = 0$ close to $\{t = 1\}$ and $\phi$ satisfies (H), the integral curves issued from points $M_0$ where $t = 1$ enter $D_0$ for $x \geq x_0$.

We similarly set

$$(Iv)(x, y, t) = \int_0^t v(x, y, s) ds$$

and write the equation $P_\alpha u = f$,

$$u + IB_1 u + I_{\tilde{S}} B_2 u + I_{\tilde{S}} IB_3 u + I^2 B_4 u + I_{\tilde{S}} I^2 B_5 u = I_{\tilde{S}} I^2 f$$

for appropriate zero order operators $B_i$ in $y$ that we now determine. Applying $\partial_t^2 \tilde{S}$ to the left, we get

$$\partial_t^2 \tilde{S} u + \partial_t^2 \tilde{S} IB_1 u + \partial_t^2 B_2 u + \partial_t B_3 u + \partial_t^2 \tilde{S} I^2 B_4 u + B_5 u = f.$$

But

$$\partial_t^2 \tilde{S} I = \varepsilon \partial_t(\partial_t \tilde{s}_0) + \partial_t \tilde{S},$$
$$\partial_t^2 \tilde{S} I^2 = 2\varepsilon \partial_t^2 \tilde{s}_0 I + 2\varepsilon \partial_t \tilde{s}_0 + \tilde{S}.$$

Hence we can take

$$B_1 = A_1, B_2 = A_2, B_4 = A_4, B_3 + \varepsilon \partial_t \tilde{s}_0 B_1 = A_3, B_5 + 2\varepsilon(\partial_t \tilde{s}_0 + \partial_t^2 \tilde{s}_0 I) B_4 = A_5.$$

A standard fixed-point argument yields the smooth flat solution $u_\alpha$.

b. One can prove, exactly as we have done for $P$, higher order energy estimates for $P_\alpha$ with constants independent of $\alpha$ (see [4]). We can find some subsequence of $u_\alpha$ converging, say, weakly in $L^2$, to some solution $u$ of $Pu = f$. Since $|u_\alpha|_s$ is bounded for all $s$, $u$ is smooth and has zero traces. It is important here that, according to Proposition IV.2.2, we do not have to decrease $\varepsilon$ with $s$. □

## V. Going back to the solution

1. In parts III and IV, we obtained a solution $(v, \phi)$ of (3.1.16) in a domain

$$D = \{-A_0 \leq s \leq M, \ \omega \in S^1, \tau_0 \leq \tau \leq \bar{\tau}_\varepsilon\}.$$

If $0 \leq \varepsilon \leq \varepsilon_k$, this solution is, say, of class $C^k$. We define now a function $\tilde{G}$ supported in $\{s \leq M\}$ by

$$\bar{\mathcal{L}}_1 \tilde{G} = \phi_s v, \ \tilde{G}(s, \omega, \tau_0) = G(s, \omega, \tau_0).$$

Note that

$$\partial_s \tilde{G} = \phi_s \tilde{v}, \ \tilde{v} = -(1 + \frac{\varepsilon^2 \phi_\tau}{2\tau})^{-1}(v - \frac{\varepsilon^2}{2\tau} \tilde{G}_\tau).$$

Thanks to condition (H), we see that (3.1.1) yields continuous functions $G$ and $F$ defined in

$$\Phi(D) = \{\phi(-A_0, \omega, \tau) \leq \sigma \leq M, \ \omega \in S^1, \ \tau_0 \leq \tau \leq \bar{\tau}_\varepsilon\}.$$

Moreover, $F$ satisfies there equation (2.2.3) $\tilde{P}(F) = 0$. Of course, $F$ and $G$ are of class $C^k$ for $\tau < \bar{\tau}_\varepsilon$, and $\nabla F$ becomes infinite at the point $\Phi(\tilde{M}_\varepsilon)$ of the upper boundary of $\Phi(D)$. Returning to the original variables and defining $\tilde{u} = \frac{\varepsilon}{r^{1/2}} G$, we obtain $\tilde{u}$ and $w$ in

$$\Omega_1 = \{\phi(-A_0, \omega, \varepsilon t^{1/2}) \leq r - t \leq M, \ \tau_0 \leq \varepsilon t^{1/2} \leq \bar{\tau}_\varepsilon\}.$$

Moreover, $w$ satisfies $P(w) = 0$ there, with $\nabla w$ becoming infinite at a point $M_\varepsilon$. Finally, for $C$ to be determined, define

$$\Omega = \{-A_0 + C(\varepsilon^2 t - \tau_0^2) \leq r - t \leq M, \tau_0 \leq \varepsilon t^{1/2} \leq \bar{\tau}_\varepsilon\} \subset \Omega_1.$$

We take $C$ big enough to have
 i) $t - A_0 + C(\varepsilon^2 t - \tau_0^2) - (t + \phi(-A_0, \omega, \varepsilon t^{1/2}))$ an increasing function,
 ii) $\Omega$ an influence domain for the operator $\partial_t^2 - \Delta - w \partial_t^2$.

In this process, we can arrange for $\Omega$ to contain $M_\varepsilon$, decreasing $\varepsilon$ if necessary (because $\phi$ is arbitrarily close to the solution for $\varepsilon = 0$, which we can arrange). The function $\tilde{u}$ satisfies the equations

$$\partial_t \tilde{u} = w, \ (\partial_t^2 - \Delta)\tilde{u} - (\partial_t \tilde{u})(\partial_t^2 \tilde{u}) = 0$$

in $\Omega$ and has the same traces as $u$ on $\{\varepsilon t^{1/2} = \tau_0\}$. By uniqueness, $u = \tilde{u}$ in $\Omega$.

The differentiability properties of $u$ close to $M_\varepsilon$ indicated in Theorems 1.1.1 and 1.1.2 follow immediately.

2. *The function $u$ does not blow up anywhere else.* The proof of the upper bound of the lifespan (1.1.1) and of Theorem 1.1.2 is already complete. It remains to prove the additionnal statement of Theorem 1.1.1 that $u$ blows up only at $M_\varepsilon$.

We recall first the *interior* asymptotic results of [1] (see Theorem 2.4 and Section 2.5). We define $S(\sigma, \omega, \tau)$ to be the solution of

$$\partial_\tau S - \frac{1}{2}(\partial_\sigma S)^2 = 0, \ S(\sigma, \omega, 0) = R^{(1)}(\sigma, \omega) + \varepsilon R^{(2)}(\sigma, \omega).$$

We then set

$$\bar{u}_a = \chi(\frac{r}{1+t})\frac{\varepsilon}{r^{1/2}} S(r - t, \omega, \tau), (\partial_t^2 - \Delta)\bar{u}_a - (\partial_t \bar{u}_a)(\partial_t^2 \bar{u}_a) = \bar{J}_a,$$

where $\chi(s)$ is zero for $s \leq 1/2$ and is one for $s \geq 3/2$. It turns out that this approximate solution satisfies

(5.2.1) $$|\partial_{x,t,\omega}^\alpha (u - \bar{u}_a)(., t)|_0 \leq C_\alpha \varepsilon^{23/9} |ln \varepsilon|$$

for $\tau$ close to $\tau_0$ and

(5.2.2) $\qquad |\partial^\alpha_{x,t,\omega} \bar{J}_a(.,t)|_0 \leq C_\alpha \varepsilon^5 |ln\varepsilon|, \ \ |\nabla^2 \bar{u}_a| \leq C\varepsilon^2$

if the norms are taken only in the "interior" domain $r - t \leq -C_0$, and $t \leq \bar{T}_\varepsilon$.

The following lemma explains how this approximate solution can be glued together with the exact solution $u$ to yield a new approximate solution $u_a$.

LEMMA V.2. *There exists an approximate solution $u_a$ with the same properties* (5.2.1) *and* (5.2.2) *as $\bar{u}_a$, and moreover $u_a = u$ for $r - t \geq -C_0$.*

*Proof.* a. Consider $G$ in a strip

$$-C_1 \leq r - t \leq -C_2$$

not containing the blowup point. We have

$$\partial_\sigma (\partial_\tau G - \frac{1}{2}(\partial_\sigma G)^2) = \varepsilon^2 [\frac{A}{\phi_s} + B](\Phi^{-1}),$$

where $A$ and $B$ are smooth. Hence, by integration, we obtain

$$\partial_\tau G - \frac{1}{2}(\partial_\sigma G)^2 = \varepsilon^2 C(\Phi^{-1}),$$

with $C$ smooth. On the other hand, the trace on $\{\tau = \tau_0\}$ of $S$ and $G$ differs by a smooth function which is $O(\varepsilon^2 ln\varepsilon)$. Hence $G - S = O(\varepsilon^2 ln\varepsilon)$.

b. If we now set

$$u_a = \chi(r-t)\bar{u}_a + (1 - \chi(r-t))u$$

for an appropriate $\chi$, we can easily check the $L^2$ estimates of the traces and of $J_a$. $\square$

We now set $u = u_a + \dot{u}$ and estimate $\dot{u}$ using standard energy inequalities between $\tau_0$ and $\bar{\tau}_\varepsilon$. When doing so, we need only control $\nabla^2 u_a$ on the support of $\dot{u}$, that is, where it is bounded. It follows that $\nabla^2 \dot{u}$, and more generally all derivatives of $\dot{u}$, are also bounded, which completes the proof.

UNIVERSITY OF PARIS-SUD, F-91405 ORSAY CEDEX, FRANCE
*E-mail address*: Serge.Alinhac@math.u-psud.fr